\documentclass{amsart}

\newtheorem{theorem}{Theorem}[section]
\newtheorem{lemma}[theorem]{Lemma}
\newtheorem{proposition}[theorem]{Proposition}
\theoremstyle{definition}

\theoremstyle{remark}

\numberwithin{equation}{section}
\theoremstyle{corollary}
\newtheorem{corollary}[theorem]{Corollary}

\def\sv{\mathcal V}
\def\su{\mathcal U}

\begin{document}

\title{Automorphisms of Coxeter groups}

\author[Patrick Bahls]{Patrick Bahls}

\address{Department of Mathematics \\ University of Illinois at 
Urbana-Champaign \\ Urbana, IL 61801}

\email{pbahls@math.uiuc.edu}



\begin{abstract}
We compute ${\rm Aut}(W)$ for any even Coxeter group whose Coxeter diagram 
is connected, contains no edges labeled 2, and cannot be separated into 
more than 2 connected components by removing a single vertex.  The 
description is given explicitly in terms of the given presentation for the 
Coxeter group and admits an easy characterization of those groups $W$ for 
which ${\rm Out}(W)$ is finite.
\end{abstract}

\thanks{The author was supported by an NSF VIGRE postdoctoral grant.}

\maketitle

\section{Introduction}

A {\it Coxeter system} is a pair $(W,S)$ consisting of a group $W$ with a 
distinguished generating set $S=\{s_i\}_{i \in I}$ for which there is a 
presentation of the form $\langle S | R \rangle$ where

$$R = \{ (s_is_j)^{m_{ij}} \ | \ m_{ij} \in \{1,2,...,\infty\}, m_{ij}=1 
\Leftrightarrow i=j, \ {\rm and} \ m_{ij}=m_{ji} \}.$$

We call $W$ a {\it Coxeter group} if it possesses a generating set $S$ for 
which $(W,S)$ is a Coxeter system.  Such a set is called a {\it 
fundamental generating set} for $W$.  If one may choose such a set $S$ so 
that all of the exponents $m_{ij}$ (for $i \neq j$) are either even or 
infinite, $W$ is said to be {\it even}, and $(W,S)$ is an {\it even 
system}.  The main result of \cite{Ba1} shows that for a given Coxeter 
group $W$, this even system is essentially unique.

We now introduce a graph-theoretical representation of a Coxeter system 
$(W,S)$.  The {\it Coxeter diagram} $\sv$ corresponding to the system 
$(W,S)$ is an edge-labeled graph whose vertex set is in one-to-one 
correspondence with the set $S$ and for which there is an edge labeled 
$m_{ij}$ between vertices $s_1$ and $s_2$ if and only if $m_{ij} < 
\infty$.  It is clear that the diagram $\sv$ is completely determined by 
the system $(W,S)$ and vice versa.  We will frequently let $\sv$ denote 
the vertex set of the Coxeter diagram; this abuse should not cause 
confusion in context.

We fix the above notation throughout the remainder of this paper, so that 
$\sv$ will always refer to the diagram corresponding to the system 
$(W,S)$.

In this paper we will concern ourselves first with Coxeter groups whose 
diagrams are connected.  (Note that if $\sv$ is not connected, then $W$ 
can be represented as a free product of the groups generated by the 
individual components of the diagram $\sv$.  More will be said about this 
in Section 6.)

We call the Coxeter group $W$ {\it rigid} if given any two systems 
$(W,S_1)$ and $(W,S_2)$, there is an automorphism $\alpha \in {\rm 
Aut}(W)$ such that $\alpha(S_1)=S_2$.  If such an automorphism can always 
be chosen from ${\rm Inn}(W)$, the group of inner automorphisms of $W$, 
then $W$ is called {\it strongly rigid}.  (In both of these cases, any two 
diagrams $\sv_1$ and $\sv_2$ for $W$ are isomorphic as edge-labeled 
graphs.)  In case $W$ is strongly rigid, the group ${\rm Aut}(W)$ has a 
very simple structure (see \cite{ChDa}):

$${\rm Aut}(W) \cong {\rm Inn}(W) \times {\rm Diag}(W)$$

\noindent where ${\rm Diag}(W)$ consists of the {\it diagram 
automorphisms} of $W$, those which are induced (in the obvious fashion) by 
symmetries of the unique diagram $\sv$ corresponding to $W$.

The goal of this paper is to describe the automorphism group ${\rm 
Aut}(W)$ of a given even Coxeter group $W$ which satisfies weaker 
conditions than strong rigidity.  This description will admit a 
finite presentation and will naturally generalize the description given 
above.

Given a Coxeter system $(W,S)$, any element of the form $wsw^{-1}$ where 
$w \in W$ and $s \in S$ is called a {\it reflection} of the system 
$(W,S)$.  (This terminology stems from the geometric action of such 
Coxeter group elements as reflections across hyperplanes in some linear 
space.)  If, for the group $W$, any two systems $(W,S_1)$ and $(W,S_2)$ 
yield the same set of reflections, we call $W$ {\it reflection 
independent}.  It is clear that strong rigidity implies reflection 
independence.  More interestingly, if $W$ is even and reflection 
independent, then it is rigid (see \cite{BaMi1}).

We say that a Coxeter system is of {\it large type} if the corresponding 
diagram has no edges labeled 2.  Using the main theorem of \cite{BaMi1} we 
conclude that any large-type even Coxeter group is reflection independent, 
and therefore rigid.  These are the groups with which we shall be 
concerned.

In \cite{Ba2}, necessary conditions for an even Coxeter group to be 
strongly rigid were given, and these conditions were shown to be 
sufficient provided that $W$ is reflection independent and that 
either $\sv$ has no simple circuits of length less than 5 (i.e.,
no ``triangles'' and no ``squares'') or $W$ is of large-type.

We will mimic the method of proof used in that paper in order to compute 
${\rm Aut}(W)$ for large-type even groups which additionally satisfy the 
following condition:

\vskip 2mm

\noindent (NVB) \ $\sv$ contains no vertex $s$ so that $\sv \setminus 
\{s\}$ consists of more than 2 connected components.

\vskip 2mm

If $\sv$ satisfies this condition, we shall say that it has {\it no vertex 
branching}, or is NVB.  We shall also say that $W$ is NVB if its diagram 
is NVB.

\begin{theorem}
Let $W$ be an even, large-type, NVB Coxeter group with connected diagram 
$\sv$.  Then ${\rm Aut}(W)$ has the form

$$G \times {\rm Diag}(W)$$

\noindent where, up to a subgroup of finite index, $G$ is a product of 
${\rm Inn}(W)$ with certain subgroups of centralizers of edges and 
vertices of the diagram $\sv$.
\end{theorem}

The precise form of the group $G$ which appears in the formula above will 
be given in Section 5, once we have developed a bit more notation and 
terminology, although the automorphisms that we construct there can be 
compared to the notion of a Dehn twist.

A {\it Dehn twist} $\phi$ of an amalgamated product $A *_C B$ 
is an automorphism satisfying $\phi(a)=a$ for all $a \in A$ and 
$\phi(b)=cbc^{-1}$ for all $b \in B$, where $c$ is a fixed central element 
in $C$.  In \cite{RS} Rips and Sela prove that for a torsion-free 
hyperbolic group $G$ the group of automorphisms of $G$ generated by the 
inner automorphisms and the Dehn twists relative to a given splitting of 
$G$ is of finite index in the overall automorphism group.  We will see 
that the group $G$ in the formula above is also generated by inner 
automorphisms and Dehn twists relative to the decomposition of $W$ as a 
certain free product with amalgamation (cf. \cite{MT} and \cite{RS}).

As a consequence of the structure of $G$ we shall also obtain

\begin{corollary}
Let $W$ be an even, large-type, NVB Coxeter group with connected diagram.  
Then ${\rm Out}(W)$ is infinite if and only if there is a vertex $s$ in 
$\sv$ such that $\sv \setminus \{s\}$ has two connected components and $s$ 
is adjacent to more than 2 other vertices in $\sv$.
\end{corollary}

In particular, if there are no cut vertices at all, we shall see that 
${\rm Aut}(W)$ has a particularly nice form:

\begin{corollary}
Let $W$ be an even, large-type, NVB Coxeter group whose connected diagram 
$\sv$ cannot be disconnected by removal of a single vertex.  Then

$${\rm Aut}(W) \cong \bigl( {\rm Inn}(W) \times Z_2^k \bigr) \times 
{\rm Diag}(W)$$

\noindent for some number $k$.
\end{corollary}

Finally, we will be able to describe the structure of ${\rm Out}(W)$, even 
in the case in which this group is infinite:

\begin{corollary}
Let $W$ be an even, large-type, NVB Coxeter group with connected diagram.  
Then ${\rm Out}(W)$ contains a subgroup of finite index isomorphic to a 
direct product of free powers of the cyclic group $Z_2$ of order 2.
\end{corollary}

Before proceeding, let us remark that a number of papers have been written 
regarding automorphisms of Coxeter groups.  (See \cite{ChDa}, \cite{FH1}, 
\cite{FH2}, \cite{HRT}, \cite{Ja}, \cite{Mu}, \cite{Tits}.)  For instance, 
\cite{FH1} and \cite{FH2} are concerned with 3-generated Coxeter groups, 
and \cite{HRT} is concerned with Coxeter groups whose diagrams are 
complete graphs.

In this paper we deal with ``larger'' Coxeter groups.  Using \cite{MT} one 
can characterize the way in which Coxeter groups decompose as free 
products with amalgamation.  The groups with which we will work here turn 
out to be the large-type even groups which are indecomposable as free 
products but which decompose in a rather nice fashion as amalgamated 
products over finite subgroups.  (By way of comparison, if $\sv$ is a 
complete graph, $W$ cannot be decomposed non-trivially as an amalgamated 
product.)

There is hope that the results of this paper might be generalized, both to 
more general classes of Coxeter groups and to other groups which arise 
as free products (with amalgamation).  Indeed, careful application of the 
results of \cite{Ba2} should yield formulas (more complicated than those 
contained in this paper) for automorphism groups of a broader class of 
even Coxeter groups than is considered here.

\section{Centralizers and junctions}

Given a subset $J \subseteq S$, we define the group $W_J \leq W$ to be the 
subgroup of $W$ generated by the generators in $J$, subject to the 
relations of the original presentation which involve only generators from 
$J$.  It is well known (see, for example, \cite{Bo}) that $(W_J,J)$ is 
itself a Coxeter system.  Clearly this system is even if $(W,S)$ is.

A subgroup $W_J$ defined in this manner is called a {\it standard 
parabolic subgroup} of $(W,S)$.  Any conjugate $wW_Jw^{-1}$ (for $w \in 
W$) of such a subgroup is called a {\it parabolic subgroup} of $(W,S)$.  

In \cite{BaMi2} the centralizers of parabolic subgroups of an even group 
$W$ are described in terms of the presentation corresponding to the given 
system $(W,S)$.  Of particular interest to us will be the centralizers of 
the standard parabolic subgroups whose generators separate the diagram 
$\sv$ into more than one connected component.  These centralizers will 
describe the amount of ``flexibility'' that we have in creating 
automorphisms of $W$.

Until Section 6 we consider only even large-type (and therefore reflection 
independent and rigid) Coxeter groups $W$ whose diagrams are connected and 
NVB.  (Rigidity of $W$ often allows us to omit mention of the set $S$.)

For any subset $J$ of $S$, let $C(J)$ denote the centralizer of $W_J$ in 
$W$.  If $J \subseteq S$ satisfies $C(J) \neq \{1\}$ and the full subgraph 
generated by the vertices $\sv \setminus J$ has more than one connected 
component, we call $J$ a {\it junction} for the Coxeter group $W$.

From Theorem 1.1 in \cite{BaMi2} we may prove the following

\begin{proposition}
Let $W_J$ be a spherical subgroup of the even large-type Coxeter group 
$W$.  Then $C(J)$ is trivial unless either

\vskip 2mm

\noindent a. \ $J=\{s_{i_1}\}$ (then $C(J) = \langle \{s_{i_1}\} \cup 
\{ (s_is_{i_1})^{\frac{m_i}{2}-1}s_i \ | \ m_i=m_{i_1i} < \infty \} 
\rangle$),
or

\vskip 2mm

\noindent b. \ $J=\{s_{i_1},s_{i_2}\}$ and $m=m_{i_1i_2} < \infty$ 
(then $C(J) = \{ 1,(s_{i_1}s_{i_2})^{m/2} \}$).
\end{proposition}

Now if $W$ is reflection independent and $\alpha \in {\rm Aut}(W)$, then 
the image $\alpha(s)$ of any fundamental generator is conjugate to some 
generator $s'$: $\alpha(s)=w_ss'w_s^{-1}$, some $w_s \in W$.  From the 
formula for ${\rm Aut}(W)$ when $W$ is strongly rigid, given in the 
previous section, we see that $w_s$ does not depend on $s$ in this case.  
In general the relationship between $w_s$ and $w_{s'}$ for $s \neq s'$ 
will be more complicated.  In order to compute ${\rm Aut}(W)$ it suffices 
to understand this relationship.

As we shall see, the junctions in the diagram $\sv$ effectively divide 
$\sv$ into subsets $S_1,...,S_k$ of the vertex set $S$ of $\sv$ so that 
$w_s=w_{s'}$ provided $s,s' \in S_i$ for some $i$.  Moreover, the 
structure of each junction will provide information about how to obtain 
one $w_{s_i}$ from an ``adjacent'' $w_{s_j}$: the quotient of two such 
words will essentially lie in the centralizer of the separating junction.  
Therefore, beginning with a ``basepoint'' in the diagram $\sv$ and an 
arbitrary conjugating word for that basepoint, we will travel throughout 
$\sv$ and compute the possible conjugating words for each vertex of $\sv$, 
allowing for the possibility that the conjugating word may be modified 
every time we cross a junction.

\section{Generators which share a common conjugating word}

Before we concern ourselves with crossing junctions, let us first 
understand those subsets of the generating set $S$ which share a common 
conjugating word $w$.

We use terminology introduced in \cite{Ba2}.  Consider the diagram $\sv$.  
A {\it simple circuit} in $\sv$ is a closed edge path $C$: $C= \{ 
[s_1s_2],...,[s_{k-1}s_k],[s_ks_1] \}$, where $s_i \neq s_j$ for $i \neq 
j$.

We say that a circuit $C=\{[s_1s_2],...,[s_ks_1]\}$ in $\sv$ is {\it 
regular} if for any two vertices $s_i$ and $s_j$ such that $|i-j|>1$, 
$s_i$ and $s_j$ are not adjacent in $\sv$.  Regular circuits are in some 
sense ``minimal'' in length.

It is essentially shown in \cite{Ba2} that if $W$ is a reflection 
independent even Coxeter group corresponding to the two diagrams $\sv$ 
and $\sv'$ and

$$C=\{[s_1s_2],...,[s_ks_1]\}$$

\noindent is a regular circuit in $\sv$, then there exists a regular 
circuit 

$$C'=\{[s'_1s'_2],...,[s'_ks'_1]\}$$

\noindent in $\sv'$ and a single word $w \in W$ such that 
$s'_i=ws_iw^{-1}$ for all $i=1,...,k$.  (In fact, from \cite{Ba1}, there 
exists a graph isomorphism $\gamma : \sv \cong \sv'$ such that 
$\gamma(s_i) = \gamma(s'_i)$, for all $i$, so the corresponding circuits 
are those which are matched up by some graph isomorphism of the diagram 
$\sv$.)

Moreover, if there is sufficient overlap between two such circuits, we can 
conclude that these two circuits must share a single conjugating word.

Suppose that $C_1$ and $C_2$ are two regular circuits in $\sv$.  Let 
$C'_1$ and $C'_2$ be the corresponding circuits in the diagram $\sv'$, and 
let $w_1$ and $w_2$ be the words guaranteed by the previous paragraphs.  
Let $J$ denote the set of vertices which lie on both cycles.  If $J$ does 
not consist of either a single vertex or two adjacent vertices, 
Proposition 2.1 shows that $C(J)=\{1\}$.  For every vertex $s \in J$, 
$w_1sw_1^{-1}=w_2sw_2^{-1}$, so that $w_1^{-1}w_2 \in C(s)$.  Thus 
$w_1^{-1}w_2 \in C(J) = \{1\}$, so $w_1=w_2$.

Furthermore, even if $J$ consists of a single vertex or a pair of adjacent 
vertices, the arguments from \cite{Ba2} show that $C_1$ and $C_2$ must 
share a common conjugating word as long as $J$ is not a junction.

These observations allow us to define a collection $\su$ of subsets of 
$S$, each of which is maximal with respect to the property that if $s_i$ 
and $s_j$ lie in the same subset $U \in \su$, then $s_i$ and $s_j$ share a 
common conjugating word.  (Although $\su$ will be a cover of $S$, it need 
not be a partition.)

To construct a typical element of $\su$, begin with a single regular 
circuit $C$, and let $U(0)=C$.  For any other regular circuit $C'$ so that 
$C \cap C' \neq \emptyset$ is not a junction, we may place the vertices of 
$C'$ in the same subset, and define $U(1)=U(0) \cup C'$.  Once we have 
considered all regular circuits which intersect $C$, we continue by 
considering those circuits which intersect any $C'$ which had been 
appended to $U(0)$, and so forth, creating a sequence $U(0),U(1),...$ of 
subsets of $S$.  Since the number of regular circuits in a given diagram 
is finite, this process must terminate with some subset $U(n) \subseteq 
S$.

We may apply this process beginning with any point that lies on a regular 
circuit.  However, we may not obtain all of the vertices in $\sv$ in this 
manner, as there may be some vertices which do not lie on a regular 
circuit.  Let $s$ be such a vertex.  If $s$ lies on {\it any} circuit, it 
is easy to show that it lies on some regular circuit, so we may assume $s$ 
does not lie on a circuit, and thus any two vertices to which is is 
adjacent must lie in different components of the full subgraph on the 
vertex set $\sv \setminus \{s\}$.  Since we are assuming that $\sv$ is 
NVB, $s$ must be incident at most 2 edges.  The vertex $s$ is a junction 
if it has degree exceeding 1.  Moreover, as will become clear in the 
following sections, such a vertex $s$ need not share a common conjugating 
word with any other generator in $S$.  Therefore the only set in the 
collection $\su$ described above which contains this vertex is $\{s\}$.

We will call each of the subsets $U$ described above a {\it unit} of the 
diagram $\sv$.  It is clear that the method described above yields a 
unique decomposition of $\sv$ into units.

\section{Conjugating words for adjacent units}

We now investigate the relationship between conjugating words for units 
which are ``near'' one another in $\sv$.

Fix $\phi \in {\rm Aut}(W)$.  Then $\phi(S)$ is another fundamental 
generating set for $W$, yielding a diagram $\sv'$ isomorphic to $\sv$.  
Reflection independence implies that for all $s_i \in S$, there is a 
unique generator $s_{\pi(i)}$ such that $\phi(s_i)$ and $s_{\pi(i)}$ are 
conjugate to one another.  It is easy to check (by using \cite{Ba1} and 
considering relations of the group $W$), that the map $\beta$ defined by 
$\beta(s_{\pi(i)})=s_i$ is a diagram automorphism.  Therefore by replacing 
$\phi$ with $\beta \circ \phi$ we can assume that $s_i$ and $\phi(s_i)$ 
are conjugate, for all generators $s_i$.  We make this assumption 
throughout this section and the next.

Let $U_i$ and $U_j$ be units in $\su$ such that either

\vskip 2mm

\noindent a. \ $J = U_i \cap U_j \neq \emptyset$ is a junction such that 
$U_i$ and $U_j$ are in different connected components of the graph $\sv 
\setminus J$, or

\vskip 2mm

\noindent b. \ there exists a unique edge $[s_is_j]$ for $s_i \in U_i$ 
and $s_j \in U_j$, and removing this edge disconnects the diagram $\sv$.  

\vskip 2mm

In either case we say that $U_i$ and $U_j$ are {\it adjacent} to one 
another.  Let $w_i$ and $w_j$ be the conjugating words associated with 
these units.  There are six cases that we now consider:

\vskip 2mm

\noindent 1. \ $|U_i|>1$, $|U_j|>1$, and $|U_i \cap U_j|=2$;

\vskip 2mm

\noindent 2. \ $|U_i|>1$, $|U_j|>1$, and $|U_i \cap U_j|=1$;

\vskip 2mm

\noindent 3. \ $|U_i|>1$, $|U_j|>1$, and $U_i \cap U_j = \emptyset$, there 
is a unique edge $[s_is_j]$ such that $s_i \in U_i$ and $s_j \in 
U_j$, and $\{s_i\}$, $\{s_j\}$ are both single-vertex junctions;

\vskip 2mm

\noindent 4. \ $|U_i|>1$, $U_j=\{s_j\}$, and there is a unique edge 
$[s_is_j]$ such that $s_i \in U_i$;

\vskip 2mm

\noindent 5. \ $U_i=\{s_i\}$, $|U_j|>1$, and there is a unique edge 
$[s_is_j]$ such that $s_j \in U_j$;

\vskip 2mm

\noindent 6. \ $U_i=\{s_i\}$, $U_j=\{s_j\}$, and $[s_is_j]$ is an edge in 
$\sv$.

\vskip 2mm

\noindent {\bf Case 1.} \ In this case both units have at least three 
vertices (by the results of the previous section) and therefore have 
trivial centralizers (by Theorem 1.1 of \cite{BaMi2}).  Let $U_i \cap U_j 
= \{s_1,s_2\}$ be a junction separating $U_i$ and $U_j$.  Then $|C(J)|=2$, 
and either $w_j=w_i$ or $w_j=w_i(s_1s_2)^{\frac{m_{12}}{2}}$.

\vskip 2mm

\noindent {\bf Case 2.} \ Now we suppose that $U_i \cap U_j = \{s\}$.  By 
Theorem 1.1 from \cite{BaMi2} we may write $w_i^{-1}w_j$ as a product of 
$s$ and the words $(s_is)^{\frac{m_i}{2}-1}s_i$, where $(s_is)^{m_i}=1$ 
and each $s_i$ is adjacent to $s$ in $\sv$.  Because $\sv$ is NVB, we may 
partition the vertices $s_i$ adjacent to $s$ into two disjoints sets, 
according to the connected component of $\sv \setminus \{s\}$ in which 
$s_i$ lies.  Suppose that $A_i$ and $A_j$ are the two connected components 
of $\sv \setminus \{s\}$ which contain $U_i \setminus \{s\}$ and $U_j 
\setminus \{s\}$, respectively.

The following lemma is fundamental.

\begin{lemma}
Suppose that $U_i$ and $U_j$ are separated by a single-vertex junction, as 
above, and that $w_i$ and $w_j$ are the respective conjugating words for 
$U_i$ and $U_j$.  Then $w_j=w_i \epsilon_s uv$, where $\epsilon_s \in 
\{1,s\}$, $u$ is a (possibly trivial) product of the words 
$(s_is)^{\frac{m_i}{2}-1}s_i$ for $s_i \in A_i$ adjacent to $s$, and $v$ 
is a (possibly trivial) product of the words $(s_js)^{\frac{m_j}{2}-1}s_j$ 
for $s_j \in A_j$ adjacent to $s$.
\end{lemma}

Thus, the ratio $w_i^{-1}w_j$ is easily controlled.

In order to prove Lemma 4.1, we make use of the following simple 
observation, which will be needed in the next section as well.

\begin{lemma}
Suppose that $U$ and $U'$ are units of $\sv$ which both contain the 
single-vertex junction $\{s\}$ such that $(U \cup U') \setminus \{s\}$ 
lies in a single connected component of $\sv \setminus \{s\}$.  Then there 
is a finite sequence of edges $[s_1s],[s_2s],...,[s_ms]$ and corresponding 
sequence of units $U=U_1,U_2,...,U_{m+1}=U'$ so that $U_k \cap U_{k+1} = 
\{s_k,s\}$ is a junction separating $U_k$ and $U_{k+1}$ for all 
$k=1,...,m$.
\end{lemma}

\begin{proof}
There is nothing to prove if $U=U'$, so we assume this is not the case.

Pick and fix vertices $t \in U$ and $t' \in U'$ such that $[st]$ and 
$[s't]$ are edges in $\sv$.  Because $t$ and $t'$ lie in the same 
component of $\sv \setminus \{s\}$, we can choose a simple path $p$ which 
does not contain $s$ and which connects $t$ to $t'$.  We assume that $p$ 
has been chosen as the shortest such path.  Concatenating this path with 
the path $\{[t's],[st]\}$ we obtain a cycle $C$.

If $C$ is regular, then $t$ and $t'$ lie in a common unit, $U''$, which is 
separated from $U$ by $\{s,t\}$ and from $U'$ by $\{s,t'\}$, and we are 
done.

Otherwise, we can ``shorten'' $C$ to form a regular cycle (compare 
\cite{Ba2}).  Since $p$ was chosen to be the shortest path from $t$ to 
$t'$ which does not contain $s$, the only way in which $C$ can fail to be 
regular is if there is some vertex $t''$ lying on $p$ for which $[st'']$ 
is an edge in $\sv$.  In this case, we may divide $C$ into two strictly 
shorter cycles (one containing $\{s,t,t''\}$ and the other $\{s,t',t''\}$) 
and induct on the length of the paths into which $p$ has been subdivided 
to yield the desired conclusion.
\end{proof}

We now prove Lemma 4.1.

\begin{proof}
We know that $w_j=w_iw$ for some word $w \in C(\{s\}$.  Let us write $w$ 
as a product

$$w = \epsilon_s \alpha_1 \beta_1 \alpha_2 \beta_2 \cdots \alpha_m 
\beta_m$$

\noindent where each $\alpha_k$ is a product of words 
$(s_ls)^{\frac{m_l}{2}-1}s_l$, $s_l \in A_i$ and each $\beta_k$ is a 
product of words $(s_ls)^{\frac{m_l}{2}-1}s_l$, $s_l \in A_j$.  We can 
clearly assume that $\alpha_k \neq 1$ for $2 \leq k \leq m$ and that 
$\beta_k \neq 1$ for $1 \leq k \leq m-1$.

Let us denote by $\phi$ the automorphism to which the conjugating 
words $w_i$ and $w_j$ correspond (that is, in particular, 
$\phi(t)=w_itw_i^{-1}$ for all $t \in U_i$, and $\phi(t)=w_jtw_j^{-1}$ for 
all $t \in U_j$).  Let $\psi$ denote the map $\phi^{-1}$.  We know that 
there exist words $w'_i$ and $w'_j$ such that $\psi(t)=w'_it{w'_i}^{-1}$ 
for all $t \in U_i$ and $\psi(t)=w'_jt{w'_j}^{-1}$ for all $t \in U_j$.  
Furthermore we can write $w'_j=w'_iw'$, where $w' \in C(\{s\})$, so that 
$w'$ can be written as a product

$$\epsilon'_s \alpha'_1 \beta'_1 \cdots \alpha'_n \beta'_n$$

\noindent for $\epsilon'_s \in \{1,s\}$ and words $\alpha'_k$ and 
$\beta'_k$ of forms similar to those of the words $\alpha_k$ and 
$\beta_k$ above.

Consider any vertex $s_l$ adjacent to $s$.  If $s_l \in A_i$, Lemma 4.2
shows that there exists a word $\bar{\alpha}_l$ which can be written 
as a product consisting solely of letters from $A_i \cup \{s\}$ such that 

$$\psi(s_l) = w'_i\bar{\alpha}_ls_l\bar{\alpha}_l^{-1}{w'_i}^{-1}. 
\eqno(1)$$

\noindent Similarly, if $s_l \in A_j$, we are guaranteed a word 
$\bar{\beta}_l$ which can be written as a product consisting solely of 
letters of $A_j \cup \{s\}$ such that

$$\psi(s_l) = w'_j\bar{\beta}_ls_l\bar{\beta}_l^{-1}{w'_j}^{-1}. 
\eqno(2)$$

Suppose that $t$ is a vertex in $U_i$.  Then

$$t = \psi \circ \phi(t) = \psi(w_itw_i^{-1}) = 
\psi(w_i)w'_it{w'_i}^{-1}\psi(w_i)^{-1}$$

\noindent so $\psi(w_i)w'_i=1$ since $C(U_i)=\{1\}$.  Now consider 
$t \in U_j$.  Here,

$$t = \psi \circ \phi(t) = \psi(w_iwtw^{-1}w_i^{-1}) = 
\psi(w_i)\psi(w)w'_iw't{w'}^{-1}{w'_i}^{-1}\psi(w)^{-1}\psi(w_i)^{-1}$$

\noindent so $\psi(w_i)\psi(w)w'_iw'=1$ since $C(U_j)=\{1\}$.  But

$$\psi(w) = \psi(\epsilon_s \alpha_1 \beta_1 \cdots \alpha_m \beta_m) = 
w'_i \epsilon_s \tilde{\alpha}_1 w' \tilde{\beta}_1 {w'}^{-1} 
\tilde{\alpha}_2 w' \tilde{\beta}_2 {w'}^{-1} \cdots \tilde{\alpha}_m w' 
\tilde{\beta}_m {w'_j}^{-1},$$

\noindent where $\psi(\alpha_k)=w'_i \tilde{\alpha}_k {w'_i}^{-1}$ and 
$\psi(\beta_k)=w'_j \tilde{\beta}_k {w'_j}^{-1}$, and equations (1) and 
(2) guarantee that $\tilde{\alpha}_k$ can be written as a product 
consisting solely of letters in $A_i \cup \{s\}$, and $\tilde{\beta}_k$ 
can be written as a product consisting solely of letters of $A_j \cup 
\{s\}$.

\noindent Using $\psi(w_i)w'_i=\psi(w_i)\psi(w)w'_iw'=1$, we obtain

$$\epsilon_s \tilde{\alpha}_1 w' \tilde{\beta}_1 {w'}^{-1} \cdots 
\tilde{\alpha}_m w' \tilde{\beta}_m = 1. \eqno(3)$$

Now we expand (3) by writing out $w'$:

$$\epsilon_s \tilde{\alpha}_1 \bigl(\epsilon'_s \alpha'_1 \cdots \beta'_n 
\bigr) \tilde{\beta}_1 \bigl({\beta'_n}^{-1} \cdots {\alpha'_1}^{-1} 
\epsilon'_s \bigr) \cdots \bigl(\epsilon'_s \alpha'_1 \cdots \beta'_n 
\bigr) \tilde{\beta}_m = 1.$$

The occurrences of $\epsilon'_s$ can be commuted to the front and 
multiplied with $\epsilon_s$ to yield a single word in $\{1,s\}$.  What 
further reduction can be performed?  Assuming $w'$ has been written in 
reduced form as a product of the words $\alpha'_k$ and $\beta'_k$, the 
only cancellation that can occur is in one of the following subwords:

\vskip 2mm

a. \ $\tilde{\alpha}_1 \alpha'_1$,

b. \ ${\alpha'_1}^{-1} \tilde{\alpha}_k \alpha'_1$, $k=2,...,m$,

c. \ $\beta'_n \tilde{\beta_k} {\beta'_n}^{-1}$, $k=1,...,m-1$, or

d. \ $\beta'_n \tilde{\beta_m}$.

\vskip 2mm

Consider the second case.  If ${\alpha'_1}^{-1} \tilde{\alpha}_k 
\alpha'_1 = 1$, then $\tilde{\alpha}_k = 1$, so

$$\alpha_k = \phi(w'_i) \tilde{\alpha}_k \phi(w'_i)^{-1} = 1.$$

But we have assumed that for $2 \leq k \leq m$, $\alpha_k \neq 1$.  Thus 
the word in the second case above is nontrivial.  Similarly we may show 
that the word in the third case above is nontrivial.

Therefore if $m>1$, what remains after cancellation is a product in words 
$(s_ls)^{\frac{m_l}{2}}s_l$ which alternates between blocks such words for 
$s_l \in A_i$ and blocks of such words for $s_l \in A_j$, and which 
represents the trivial element.  Because the groups that we are 
considering satisfy the $C'(\frac{1}{6})$ small cancellation condition, 
any non-trivial word representing the trivial element must contain more 
than half of a relator appearing in the symmetrized presentation for the 
group.  (See \cite{LS} for more details.)  However, this is clearly not 
the case if $m>1$.  Therefore $m=1$, and the word $w$ can be written 
$\epsilon_s uv$ for words $u$ and $v$ described in the statement of the 
lemma.

(Note that we have assumed that $w$ (resp. $w'$) begins with some 
word $\alpha_1$ (resp. $\alpha'_1$) and ends with some word $\beta_m$ 
(resp. $\beta'_n$); a moment's thought should convince the reader that 
the other possibilities are similar.)
\end{proof}

\noindent {\bf Case 3.} \ Here, although $U_i$ and $U_j$ no longer 
overlap, there is a {\it bridge} $[s_is_j]$ between them in $\sv$.  (That 
is, removing this edge disconnects the remaining graph.)

Because the edge $[s_is_j]$ does not lie on a cycle in $\sv$, we cannot 
directly apply the methods of \cite{Ba2}.  However, by the assumptions 
made in the first paragraph of this section, the subgroups $H = 
W_{\{s_i,s_j\}}$ and $W_{\{\phi(s_i),\phi(s_j)\}}$ are conjugate, so 
there is a word $w$ satisfying $ws_iw^{-1}=\phi(s_i)$ and 
$wxs_jx^{-1}w^{-1}=\phi(s_j)$ for some word $x$ in the letters $s_i$ and 
$s_j$ such that $s_i$ and $xs_jx^{-1}$ generate $H$.

Clearly $w_i^{-1}w \in C(\{s_i\})$ and $x^{-1}w^{-1}w_j \in C(\{s_j\})$, 
so that $w_i^{-1}w_j = \hat{s}_ix\hat{s}_j$ for some $\hat{s}_i \in 
C(\{s_i\})$ and $\hat{s}_j \in C(\{s_j\})$.  We must understand $x$.  Let 
$n$ be the label on the edge $[s_is_j]$.

Because $s_i$ and $xs_jx^{-1}$ generate $H$, $s_ixs_jx^{-1}=(s_is_j)^k$ 
for some $k$ such that $(k,n)=1$ (where $(k,n)$ denotes the greatest 
common divisor of $k$ and $n$).  For a given $k \in \{1,...,n-1\}$ such 
that $(k,n)=1$, there are four geodesic (that is, shortest among words 
representing the same group element) words $x$ which satisfy 
$s_ixs_jx^{-1} = (s_is_j)^k$.  Namely, for a given $k$, $x$ lies in the 
set

$$\{ x_1 = (s_js_i)^{\frac{k-1}{2}}, x_2 = 
(s_js_i)^{\frac{k-1}{2}}s_j, x_3 = (s_is_j)^{\frac{n-k-1}{2}}s_i, x_4 = 
(s_is_j)^{\frac{n-k+1}{2}} \}. \eqno(4)$$

Conjugating $s_j$ by any one of these elements $x_l(k)$ yields the same 
element of $W$.  However, since for the given automorphism $\phi$, 
$\phi(s) = w_jsw_j^{-1}$ must hold for every $s \in U_j$, it is clear that 
different choices of $x_l(k)$ for the same $k$ will yield different maps 
$\phi$.

Suppose that a given automorphism $\phi_l$ satisfies

$$\phi(s) = w_isw_i^{-1}$$

\noindent for $s \in U_i$ and

$$\phi(s) = w_iux_l(k)vsv^{-1}x_l(k)^{-1}v^{-1}w_i^{-1}$$

\noindent for $s \in U_j$, for some $k,l$, and for some words $u$ and $v$ 
as described in Case 2.  (That is, $u$ is a product of words 
$(ss_i)^{\frac{m}{2}-1}s$ for $s \in U_i$ and $v$ is a product of words 
$(ss_j)^{\frac{m}{2}-1}s$ for $s \in U_j$.)  What can $k$ and $l$ be?  The 
following lemmas effectively guarantee that for any choice of $k$ and $l$ 
as above, there are numbers $k'$ and $l'$ such that an automorphism 
corresponding to these choices is inverse to the first.  The proof of the 
first involves elementary number theory and is omitted.  The proof of the 
second requires tedious but straightforward computation and will also be 
omitted.

\begin{lemma}
Let $k$ and $n$ be as above.  Then there is a unique number $k'$, $1 \leq 
k' \leq n-1$, such that $(k',n)=1$ and $kk'+1=dn$ for some integer $d$.
\end{lemma}

\begin{lemma}
Let $\gamma_l$, $u$, $v$, $k$, $k'$, and $d$ be as above.  Let $u'=u^{-1}$ 
and $v'=v^{-1}$.  Let $\gamma'_{l'}$ be any map of $W$ satisfying

\vskip 3mm

1. \ $\gamma'_{l'}(s)={w'}_is{w'}_i^{-1}$ for $s \in U_i$,

2. \ $\gamma'_{l'}(s)={w'}_iu'x_{l'}(k')v'svx_{l'}(k')^{-1}u{w'}_i^{-1}$ 
for $s \in U_j$, and

3. \ $\gamma'_{l'}(w_i)w'_i=1$.

\vskip 3mm

\noindent If $d$ is odd, then

$$\gamma_1 \circ \gamma'_4(s) = \gamma'_4 \circ \gamma_1(s) = \gamma_2 
\circ \gamma'_3(s) = \gamma'_3 \circ \gamma_2(s) = s$$

\noindent for all $s \in U_i \cup U_j$.  If $d$ is even, then

$$\gamma_l \circ \gamma'_l(s) = \gamma'_l \circ \gamma_l(s) = s$$

\noindent for all $s \in U_i \cup U_j$ and $l=1,2,3,4$.
\end{lemma}

Therefore, any choice of $x_l(k)$ yields a map which is ``invertible'' 
when restricted to $U_i \cup U_j$.  We have $w_j=w_iux_l(k)v$, for any $k$ 
such that $(k,n)=1$, any $l \in \{1,2,3,4\}$, and any $u$ and $v$ as 
above.  (Note that we include neither $\epsilon_{s_i} \in \{1,s_i\}$ nor 
$\epsilon_{s_j} \in \{1,s_j\}$ as was done before, as $s_ix_l(k) = 
x_{l'}(n-k)$ for some $l'$, and $x_l(k)s_j = x_{l'}(k)$ for some $l'$.)

\vskip 2mm

\noindent {\bf Case 4.} \ Now $|U_i|>1$ and $U_j=\{s_j\}$.  As was argued 
in Section 3, there exists a unique vertex $s_i \in U_i$ such that 
$[s_is_j]$ is an edge in $\sv$, and this vertex $s_i$ is a junction 
separating $U_i$ and $U_j$.  We can now show, arguing as before, that 
$w_j=w_iux_1(k)$, where $u$ is as before, and we select the single 
element $x_1(k)$ from each set of four which elements conjugates $s_j$ to 
$(s_js_i)^{k-1}s_j$.  (Recall that the effect of conjugating $s_j$ by any 
one of these four elements is the same, for a given $k$.)

\vskip 2mm

\noindent {\bf Case 5.} \ Now $U_i=\{s_i\}$ and $|U_j|>1$.  We can argue 
as above to show that $w_j=w_iux_l(k)v$, where $u$, $v$, and $x_l(k)$ are 
defined as before.  (Note that if $s_i$ has degree 1, $u$ must be 
trivial.)

\vskip 2mm

\noindent {\bf Case 6.} \ Finally, suppose that $|U_i|=|U_j|=1$.  In this 
case we may show that $w_j=w_iux_1(k)$, where $u$ and $x_1(k)$ are as 
before.  As in Case 5, $u=1$ if $s_i$ has degree 1.

\section{The unit graph and the structure of an automorphism}

We now understand how to modify our conjugating element when passing from 
one unit to any adjacent unit.  In this section we piece together this 
local information to obtain an arbitrary automorphism.

We define the {\it unit graph} $\Gamma=\Gamma(\sv)$ of the diagram $\sv$ 
as the (unlabeled) graph whose vertex set is $\su$ and for which there is 
an edge $[U_iU_j]$ between units $U_i$ and $U_j$ whenever $U_i$ and $U_j$ 
are adjacent.

It is easy to see that $\Gamma$ is connected.  (One can connect any two 
units by keeping track of the units that are entered in tracing a path 
from a vertex in one unit to a vertex in the other.)  In $\Gamma$, choose 
any spanning tree $T$.

Let $\{s\}$ be any single-vertex junction in $\sv$.  Because $\sv$ is NVB, 
we may divide the units of $\sv$ which contain $s$ into two sets, 
depending the connected component of $\sv \setminus \{s\}$ in which a 
given unit lies.  Denote these two subsets of $\su$ by $\su_1(s)$ and 
$\su_2(s)$.  We claim that we can modify the spanning tree $T$ to obtain a 
new spanning tree $T'$ which has the property that any time $\{s\}$ is a 
single-vertex junction, there is exactly one pair $(U_1,U_2) \in \su_1(s) 
\times \su_2(s)$ for which $[U_1U_2]$ is an edge in $T'$.

Indeed, let $s$ be as above.  Clearly $T$ contains at least one edge as 
described above to begin with (this is because $\{s\}$ is a junction).  
Denote it by $[U_1U_2]$.  By the definition of $\su_i(s)$ ($i=1,2$), it is 
clear that every element $U \in \su_1(s)$ is connected to $U_1$ in 
$\Gamma$, and similarly every element $U \in \su_2(s)$ is connected to 
$U_2$ in $\Gamma$.  Therefore if $[U'_1U'_2]$ is any other edge in $T$ for 
$(U'_1,U'_2) \in \su_1(s) \times \su_2(s)$, we can eliminate this edge 
from $T$ and maintain connectivity by adding to $T$ (if necessary) paths 
from $U'_1$ to $U_1$ and from $U'_2$ to $U_2$.  (Essentially, one simply 
creates a spanning tree for the subgraph $\su_1(s) \cup \su_2(s)$ which 
uses the prescribed edge $[U_1U_2]$.)

The reader should convince himself or herself that this modification can 
be performed independently for all single-vertex junctions.  (This is so 
because no elements of $\su_1(s)$ are separated from each other by 
single-vertex junctions, by Lemma 4.2.)

\vskip 2mm

\noindent {\bf Remark.} \ Why perform this modification?  We shall use the 
tree $T'$ to define an automorphism by stepping from one vertex of $T'$ to 
an adjacent vertex and modifying the conjugating word in the manner 
appropriate for the edge so traversed.  Suppose $s \in U_i \cap U_j$ and 
$U_i, U_j \in \su_1(s)$.  Then Lemma 4.2 implies that the conjugating 
words $w_i$ and $w_j$ must differ only by a product of words 
$(ss_k)^{\frac{m_k}{2}}$ for $s_k$ lying in the same connected component.  
In order that we yield such a modification, we must not be allowed to 
``cross over'' a single-vertex junction and then ``cross back''.

\vskip 2mm

Now pick and fix any vertex $U_0$ in $T'$ such that $U_0$ is not a 
single-vertex unit $\{s\}$ for $s$ a vertex of degree 1.  (It is clear 
this choice can always be made.)  Call $U_0$ the {\it basepoint} of $T'$.  
Since $T'$ is a tree, there is a unique geodesic (shortest) path from 
$U_0$ to any other vertex $U$ of $T'$.  We assign an orientation to the 
edges of $T'$.  If $[U_iU_j]$ is any edge in $T'$, we define $U_i$ to be 
the initial point of $[U_iU_j]$ and $U_j$ to be the terminal point of 
$[U_iU_j]$ if the geodesic path from $U_0$ to $U_j$ passes through $U_i$.

We now label the (oriented) edges of $T'$ with elements of $W$.  For an 
edge $[U_iU_j]$, we will denote its label by $\phi_{[U_iU_j]}$.

The edges of $T'$ are of six types, corresponding exactly to the six cases 
considered in the previous section.  For any junction $J$ separating $U_i 
\setminus J$ and $U_j \setminus J$, let $A_i(J)$ be the connected 
component of $\sv \setminus J$ containing $U_i \setminus J$, and let 
$A_j(J)$ be the connected component of $\sv \setminus J$ containing $U_j 
\setminus J$.

\vskip 2mm

\noindent 1. \ If $|U_i|>1$, $|U_j|>1$, and $U_i \cap U_j = \{s_1,s_2\}$ 
is a junction separating $U_i$ and $U_j$, then $\phi_{[U_iU_j]}=u$, where 
$u \in \{1,(s_1s_2)^{\frac{m_{12}}{2}}\}$.

\vskip 2mm

\noindent 2. \ If $|U_i|>1$, $|U_j|>1$, and $U_i \cap U_j = \{s\}$ is a 
junction separating $U_i$ and $U_j$, then $\phi_{[U_iU_j]}=\epsilon_suv$, 
where $\epsilon_s \in \{1,s\}$, $u$ is a product of words 
$(s_is)^{\frac{m_i}{2}-1}s_i$ for $s_i$ in $A_i(\{s\})$, and $v$ is a 
product of words $(s_js)^{\frac{m_j}{2}-1}s_j$ for $s_i$ in $A_j(\{s\})$.

\vskip 2mm

\noindent 3.\ If $|U_i|>1$, $|U_j|>1$, $U_i \cap U_j = \emptyset$, and 
$[s_is_j]$ is an edge labeled $n$ for $s_i \in U_i$ and $s_j \in U_j$, 
then $\phi_{[U_iU_j]} = ux_l(k)v$, where $u$ is a product of words 
$(ss_i)^{\frac{m}{2}-1}s$ for $s \in A_i(\{s_i\})$, $v$ is a product of 
words $(ss_j)^{\frac{m}{2}-1}s$ for $s \in A_j(\{s_j\})$, and $x_l(k)$ is 
one of the words in (4) from Section 4, for some $k$ such that $(k,n)=1$ 
and some $l=1,2,3,4$.

\vskip 2mm

\noindent 4. \ If $|U_i|>1$, $U_j=\{s_j\}$, and $s_i \in U_i$ is such that 
$[s_is_j]$ is an edge in $\sv$ labeled $n$, then 

\vskip 2mm

a. \ if $s_j$ has degree 1, $\phi_{[U_iU_j]}=ux_1(k)$, and 

b. \ if $s_j$ has degree 2, $\phi_{[U_iU_j]}=ux_l(k)$,

\vskip 2mm

\noindent where $u$ is a product of words $(ss_i)^{\frac{m}{2}-1}s$ for $s 
\in A_i(\{s_i\})$, and $x_l(k)$ is as in (4) from Section 4, for some $k$ 
such that $(k,n)=1$, $l=1,2,3,4$.

\vskip 2mm

\noindent 5. \ If $U_i=\{s_i\}$, $|U_j|>1$, and $s_j \in U_j$ is such that 
$[s_is_j]$ is an edge in $\sv$ labeled $n$, then $\phi_{[U_iU_j]} = 
x_l(k)v$, where $x_l(k)$ and $v$ are as in Case 3, where we allow only 
$l=1$ and $l=2$.

\vskip 2mm

\noindent 6. \ If $U_i=\{s_i\}$ and $U_j=\{s_j\}$ and $[s_is_j]$ is an 
edge labeled $n$ in $\sv$, then

\vskip 2mm

a. \ if $s_j$ has degree 1, then $\phi_{[U_iU_j]}=x_1(k)$, and

b. \ if $s_j$ has degree 2, then $\phi_{[U_iU_j]}=x_l(k)$,

\vskip 2mm

\noindent where $x_l(k)$ is as in Case 3, for either $l=1$ or $l=2$.

\vskip 2mm

We note that the labels we give each sort of edge differ only slightly 
from the corresponding ratios $w_i^{-1}w_j$ computed in the previous 
section.  In order to explain the differences, we will need the following 
simple lemma, whose proof consists of direct computation.

\begin{lemma}
Let the word $x_l(k)$ be defined as above, relative to the edge $[s_is_j]$ 
labeled $n$.  Then

\vskip 2mm

a. \ $x_3(k)=s_i x_1(n-k)$ and $x_4(k)=s_i x_2(n-k)$ for all $k$ such that 
$(k,n)=1$, and

\vskip 2mm

b. \ for all $k$ such that $(k,n)=1$, $l \in \{1,2,3,4\}$, and $z \in \{1, 
s_j, (s_is_j)^{\frac{n}{2}-1}s_i, (s_is_j)^{\frac{n}{2}} \}$, there exists 
$l' \in \{1,2,3,4\}$ such that $x_l(k)z=x_{l'}(k)$.
\end{lemma}

In Case 5, we seem to have lost generality by removing the term $u$ and by 
allowing only $l=1$ and $l=2$.

However, if $s_i$ is a vertex with degree 1, $u=1$ must hold.  When $s_i$ 
has degree 2, $u \in \{ 1, s_i, (ss_i)^{\frac{m}{2}-1}s, 
(ss_i)^{\frac{m}{2}} \}$ (where $[ss_i]$ is an edge labeled $m$, $s \neq 
s_j$).  Further, Lemma 5.1 allows us to factor a single letter $s_i$ out 
of either $x_3(k)$ or $x_4(k)$.  This letter can then be multiplied with 
$u$; the resulting product may then be absorbed by the term $x_l(k)$ which 
occurs in the word $\phi_{[\{s\}\{s_i\}]}$.  Lemma 5.1 now guarantees that 
such absorption does not alter the element to which $s_i$ is conjugated.  
(Also note that it is clear that the edge $[\{s\}\{s_i\}]$ must appear in 
$T'$.)  

In Case 6, we may omit the term $u$ for the same reason.  In this case 
and in Case 4, we must also be ready to absorb a term which may come 
from the following edge in $T'$.  We allow $l=1$ and $l=2$ in Case 6 to 
account for this absorption.  If either $x_3(k)$ or $x_4(k)$ in {\it this} 
term arises as a result of absorption, we may factor out a single letter 
$s_i$ and push it into the {\it previous} edge.  This process, analogous 
to ``carrying'' in arithmetic, must terminate when we reach a junction as 
in Case 4, where we allow any value of $l$ in $\{1,2,3,4\}$.

\vskip 2mm

\noindent {\bf Remark.} \ These precautions are taken in order that the 
description of an automorphism given below is unique.

\vskip 2mm

Finally, we label the vertex $U_0$ with any element $\phi_0$ of $W$.

Every such labeling defines an automorphism $\phi$ in the following 
fashion.  Fix a labeling as described above and let $s$ be a generator in 
$S$, and choose any unit $U$ which contains $s$.  Let 
$\{[U_0U_1],[U_1U_2],...,[U_{l-1}U_l]\}$ be the unique geodesic in $T'$ 
from $U_0$ to $U_l=U$.  Then define

$$\phi_s = \phi_0 \prod_{i=0}^{l-1} \phi_{[U_iU_{i+1}]}$$

\noindent and

$$\phi(s) = \phi_s s \phi_s^{-1}.$$

First note that the choice of the unit $U$ containing $s$ is not 
important.  Indeed, suppose $s \in U_i \cap U_j$.  Then $s$ is contained 
in a junction separating those two units and one sees easily that 
$\phi_{[U_iU_j]}$ must commute with $s$.

By considering group relations, it is easy to see that this map is a 
homomorphism.  To see that $\phi$ is an automorphism, we now indicate 
a formula for the composition of two such maps, and a formula for the 
inverse of a given map.

Let $\phi$ and $\phi'$ be two homomorphisms defined as above.  Let 
$\psi=\phi' \circ \phi$.  We compute $\psi$ by computing $\psi_{[U_iU_j]}$ 
for each edge $[U_iU_j]$, depending on the type of each edge:

\vskip 2mm

\noindent 1. \ $\psi_{[U_iU_j]} = \phi'_{[U_iU_j]}\phi_{[U_iU_j]}$.

\vskip 2mm

\noindent 2. \ If $\phi_{[U_iU_j]} = \epsilon_suv$ and $\phi'_{U_iU_j]} = 
\epsilon'_su'v'$, then $\psi_{[U_iU_j]} = \epsilon_s\epsilon'_suu'v'v$.

\vskip 2mm

\noindent 3. \ If $\phi_{[U_iU_j]} = ux_l(k)v$ and $\phi'_{[U_iU_j]} = 
u'x_{l'}(k')v'$, then $\psi_{[U_iU_j]} = uu'x_{l \circ l'}(kk')v'v$, where 
$\{1,2,3,4\}$ is made a group isomorphic to the 4-element group $V$ with 
$1$ as its unit and product $\circ$.

\vskip 2mm

\noindent 4. \ If $\phi_{[U_iU_j]} = ux_l(k)$ and $\phi'_{[U_iU_j]} = 
u'x_{l'}(k')$, then $\psi_{[U_iU_j]} = uu'x_{l \circ l'}(kk')$, with $l 
\circ l'$ defined as in Case 3.

\vskip 2mm

\noindent 5. \ If $\phi_{[U_iU_j]} = x_l(k)v$ and $\phi'_{[U_iU_j]} = 
x_{l'}(k')v'$, then $\psi_{[U_iU_j]} = x_{l \circ l'}(kk')v'v$, with $l 
\circ l'$ defined as in Case 3.

\vskip 2mm

\noindent 6. \ If $\phi_{[U_iU_j]} = x_l(k)$ and $\phi'_{[U_iU_j]} = 
x_{l'}(k')$, then $\psi_{[U_iU_j]} = x_{l \circ l'}(kk')$, with $l \circ 
l'$ defined as in Case 3.

\vskip 2mm

Finally, given $\phi_0$ and $\phi'_0$ corresponding to $\phi$ and 
$\phi'$, we compute $\psi_0=\phi'(\phi_0)\phi'_0$.  This completes the 
description of the composition $\psi$.

Using Lemma 4.4 in order to compute ``inverses'' for the terms $x_l(k)$, 
and letting $u'=u^{-1}$ and $v'=v^{-1}$ throughout, it is now easy to 
compute the inverse $\phi^{-1}$ of a given map $\phi$  (The most difficult 
part is to ``invert'' $\phi_0$.)  Therefore, every map so constructed is 
an automorphism of $W$.  

It is clear by the arguments given in Section 4 that all automorphisms 
$\phi$ which satisfy $\phi(s) = w_ssw_s^{-1}$ for all $s \in S$ (where 
$w_s$ depends on $s$) can be obtained in this fashion.  Moreover it is not 
difficult to show that if $\phi$ and $\phi'$ are defined as above, then 
$\phi=\phi'$ if and only if $\phi_0=\phi'_0$ and $\phi_{[U_iU_j]} = 
\phi'_{[U_iU_j]}$ for all edges $[U_iU_j]$ in $T'$.

Given any automorphism of $W$, we obtained an automorphism in the 
collection $G$ of automorphisms described above by composing it with an 
element of ${\rm Diag}(W)$.  The intersection $G \cap {\rm Diag}(W)$ is 
trivial, so that every element of ${\rm Aut}(W)$ can be written uniquely 
as such a composition.  Thus ${\rm Aut}(W)$ is a semidirect product of $G$ 
by ${\rm Diag}(W)$, and we have proved the main theorem.

It is easy to see that the inner automorphisms and diagram automorphisms, 
along with automorphisms for which $\phi_{[U_iU_j]}=1$ for all but a 
single edge $[U_iU_j]$ in $T'$, generate ${\rm Aut}(W)$.  Thus ${\rm 
Aut}(W)$ is finitely generated.  It is not difficult to explicitly 
describe a finite presentation for ${\rm Aut}(W)$ in terms of these 
generators, given the formula for composition given above.

The corollaries given in the introduction are easily proven.  Because 
$|{\rm Diag}(W)| < \infty$, the only way in which ${\rm Out}(W)$ can be 
infinite is if there are infinitely many choices for $\phi_{[U_iU_j]}$, 
for some edge $[U_iU_j]$ in $T'$.  This is clearly the case only when 
there is a cut vertex which is adjacent to more than two 2 vertices in 
$\sv$.  Thus Corollary 1.2 follows.

In fact, consider the subgroup $F$ of ${\rm Aut}(W)$ generated by those 
$\phi$ for which $\phi_{[U_iU_j]}=1$ for edges of types 1 and 6, and 
$\epsilon_s=1$ and $x_l(k)=1$ for all occurrences of these words in edges 
of types 2, 3, 4, and 5.  (The group $F$ has finite index in ${\rm 
Out}(W)$.)  Then $F$ is clearly generated by a set of involutions whose 
cardinality can be read immediately from $\sv$ by counting the number of 
vertices adjacent to cut vertices.  It is not difficult to see that $F$ is 
a direct product of free powers of $Z_2$, proving Corollary 1.4.

It is also easy to compute the order of ${\rm Out}(W)$ if it this 
group is finite.  If there is a cut vertex, Corollary 1.2 and the fact 
that $\sv$ is NVB imply that $\sv$ must be a path of length $r$ which does 
not intersect itself.  Let $n_i$ be the label on the $i$th edge, 
proceeding from one of the endpoints of this path.  Then $|{\rm Out}(W)| = 
\delta n$, where

$$n=2^{r-2} \prod_{i=1}^r \phi(n_i)$$

\noindent ($\phi$ is Euler's totient function), and $\delta \in \{1,2\}$ 
(depending on whether the path $\sv$ has one or two symmetries).

If there are no cut vertices, then the only junctions are edges.  If $k$ 
is the number of units into which $\sv$ is divided by these edges, it is 
easy to see that ${\rm Aut}(W)$ has the structure given in Corollary 1.3, 
and $|{\rm Out}(W)| = 2^k |{\rm Diag}(W)|$.

\section{Free products}

As promised in the introduction, let us briefly consider the issue of 
free products of Coxeter groups.  Any Coxeter group can easily be written 
as a free product of freely indecomposable Coxeter groups.  Indeed, a 
Coxeter group is freely indecomposable if and only if its diagram is 
connected.  Moreover, given two Coxeter groups $W_1$ and $W_2$ with 
diagrams $\sv_1$ and $\sv_2$, respectively, a diagram for $W_1 * W_2$ is 
obtained by taking the the disjoint union of $\sv_1$ and $\sv_2$.

First let us consider the case of a Coxeter group $W$ which decomposes as 
free product of exactly two freely indecomposable Coxeter groups (thus its 
diagram $\sv$ has exactly two connected components).  We write $W \cong 
W_1 * W_2$ for some Coxeter groups $W_i$ with diagrams $\sv_i \subseteq 
\sv$, $i=1,2$ whose disjoint union is $\sv$.

Let us assume that $W$ is a large-type even rigid Coxeter group, and that 
both $W_1$ and $W_2$ are strongly rigid.  (It is easy to see that $W$ 
itself cannot be strongly rigid.)  Consider any automorphism $\phi$ of $W$ 
such that for every $s$ there is an element $w \in W$ such that 
$\phi(s)=wsw^{-1}$.  (As in Section 4, any $\phi \in {\rm Aut}(W)$ can be 
composed with a diagram automorphism to yield such an automorphism.)

By strong rigidity there are elements $w_i \in W_i$ such that 
$\phi(s_i)=w_is_iw_i^{-1}$ for $s_i \in W_i \cap S$, $i=1,2$.  In order to 
completely describe $\phi$, we need only compute the ratio $w_1^{-1}w_2$.  
An argument almost identical to the proof of Lemma 4.1 yields the 
following result.

\begin{lemma}
Let $W$, $W_i$, and $w_i$ be as above.  Then $w_1^{-1}w_2=u_1u_2$, where 
$u_i \in W_i$ for $i=1,2$.
\end{lemma}

Therefore every automorphism which takes a generator $s$ to a conjugate of 
$s$ is described by a triple $(w,u_1,u_2)$.  If $\phi$ is given by 
$(w,u_1,u_2)$ and $\phi'$ by $(w',u'_1,u'_2)$, it is easy to show that 
composition is given by the formula

$$(w',u'_1,u'_2) \circ (w,u_1,u_2) = (\phi'(w)w',u_1u'_1,u'_2u_2).$$

\noindent From this formula inverses can be easily computed.  If neither 
$W_1$ nor $W_2$ is finite, two distinct triples correspond to distinct 
automorphisms.  (In general, one must perform a further quotient by the 
centralizers $C(W_i)$ for $i=1,2$, but $C(W_i)$ is trivial if $W_i$ is 
infinite.)  We obtain the following result.

\begin{theorem}
Let $W$ be a large-type even rigid Coxeter group which decomposes as the 
free product of exactly two strongly rigid Coxeter groups $W_1$ and 
$W_2$.  Then

$${\rm Aut}(W) \cong G \times {\rm Diag}(W)$$

\noindent where $G$ is (element-wise) a product of $W$, $W_1$, and $W_2$, 
with multiplication given above.  Also, ${\rm Out}(W)$ is finite if and 
only if both $W_1$ and $W_2$ are finite; that is, both are isomorphic 
either to $Z_2$ or to a strongly rigid dihedral group $D_n$, for some $n$.
\end{theorem}

The statement regarding finiteness of ${\rm Out}(W)$ is clear.

The preceding arguments can also be used in a much more general setting in 
order to establish the following result.

\begin{theorem}
Let $G \cong G_1 * G_1$ where ${\rm Out}(G_i)$ is finite for $i=1,2$.  
Then ${\rm Aut}(G)$ contains a subgroup of finite index which is 
isomorphic to a product of $G$, $G_1$, and $G_2$, with composition given 
by

$$\phi' \circ \phi = (g',g'_1,g'_2) \circ (g,g_1,g_2) = 
(\phi'(g)g',g_1g'_1,g'_2g_2)$$

\noindent for $g \in G$ and $g_i \in G_i$, $i=1,2$.
\end{theorem}

To finish, let us now consider a Coxeter group $W$ which decomposes as a 
free product of finitely many freely indecomposable Coxeter groups: $W = 
W_1 * \cdots * W_k$.  Let us furthermore assume that each $W_i$ 
($i=1,...,k$) is even, of large type, and is NVB.  Each $W_i$ is finitely 
presented, and the results from Section 5 guarantee that ${\rm Aut}(W)$ is 
finite presented for every $i$.  Therefore, we appeal to \cite{Gil}, which 
produces an explicit (and finite!) presentation for the group ${\rm 
Aut}(W)$.  We obtain the following

\begin{theorem}
Let $W$ be the free product of the even, large-type, NVB Coxeter 
groups $W_i$.  Then ${\rm Aut}(W)$ is finitely presented.  (Moreover, one 
may explicitly compute a presentation.)
\end{theorem}

If $k \geq 3$, ${\rm Out}(W)$ will always be infinite.

\end{document}